\let\documentclass\relax\documentclass 
\input AHTOHFIE.STY
\hfuzz4.8pt

\UDC{%
515.142.321
+512.543.16
+519.179.1
%
%
}


\title{%
%
Small non-Leighton two-complexes
}
\author{%
Natalia S. Dergacheva
\quad
Anton A. Klyachko
}
\address{%
Faculty of Mechanics and Mathematics of Moscow State University
\\
Moscow 119991, Leninskie gory, MSU.
\\
Moscow Center for Fundamental and Applied Mathematics
\\
nataliya.dergacheva@gmail.com
\quad
klyachko@mech.math.msu.su
}

\grants{\RSF 22-11-00075}

\abstract{%
How many 2-cells must two finite
CW-complexes have to admit a common, but not finite common, covering?
Leighton's theorem says that both complexes must have
2-cells. We construct an almost (?) minimal example with two 2-cells
in each complex.
}

\s 0.
Introduction

\proclaim{Leighton's theorem} {\rm[Lei82]}.
If two finite graph have a common covering, then they have a common finite
covering.

\noindent
Alternative proofs and various generalisations of this result can be
found, e.g., in
[Neu10],
[BaK90],
[SGW19],
[Woo21],
[BrS21],
and references therein.

Does a similar result hold for any CW-complexes,
i.e.
\disp{\hfuzz32pt
\sl
is it true
that, if, for finite CW-complexes $K_1$ and $K_2$, there exist a
CW-complex $K$ and
\spacing{cellular}
coverings
$K_1\ot K\to K_2$,
then there exists a finite CW-complex $K$ with this property?
}%
This natural question was
posed
(in other terms)
in
[Tuc90] and [AFS91].
Notice the cellularity requirement.
Surely, we would obtain an equivalent
question if we replace this condition with
a formally
stronger \emph{combinatorialness} one:
the image of each cell is a cell.
However, without
the cellularity condition,
the answer would be negative:
indeed,
the torus and the genus-two surface
have no finite common coverings (as the
fundamental group the genus-two orientable surface
$\pres<x,y,z,t|[x,y][z,t]=1>$ contains no abelian subgroups of finite
index), while
the universal coverings of these surfaces are
homeomorphic, because they are the plane.
The cellularity condition rules out such examples:
if we take, e.g., the standard one-vertex cell structures
on the torus and genus-two surface,
then, on the covering plane, we
obtain
\-
the usual square lattice on the (Euclidean) plane (in the torus-case)
\-
and an octagonal lattice on the (Lobachevskii) plane (in the genus-two
case)

\enditem
(i.e. though the universal coverings are homeomorphic,
the cell structure on them are principally different).
This example cannot be saved by
a complication of the cell structures on the torus and
genus-two surface
(as was noted in~[Tuc90] and [AFS91]; the authors of [AFS91] even
conjectured that the answer to the (cellular version of) the question
is positive).

Nevertheless,
the answer turned out to be negative as was shown in~[Wis07]
(and actually, much earlier in~[Wis96]);
the complexes $K_1,\;K_2$
forming such a
\emph{non-Leighton}
pair from [Wis07]
contain as few as six 2-cells each.
In [JaW09], this number was reduced to four:%
\fn{%
although the authors of [JaW09]
did not pursue this purpose;
it was rather a
byproduct of their results.}
\disp{\sl\narrower\narrower\narrower
there exist two two-complexes
containing four 2-cells each
that
have a common covering but have not finite common coverings.
}%
(Henceforth, we omit the prefix
``CW-" and word
``cellular":
a \emph{complex}
means a CW-complex,
and
all mapping between complexes
are assumed to be cellular in this paper.)
The non-Leighton complexes~$K_1$~and~$K_2$
from [JaW09]
are
the \emph{standard complexes} of the following group presentations
$\Gamma_i$,
i.e. one-vertex complexes with edges corresponding to the generators
and 2-cells attached by the relators:
$$
\Gamma_1=F_2\times F_2
=
\pres<a,b,x,y|
[a,x]
,\;
[a,y]
,\;
[b,x]
,\;
[b,y]
>
\hbox{ and }
\Gamma_2
=
\pres<a,b,x,y|
axay
,\;
ax^{-1}by^{-1}
,\;
ay^{-1}b^{-1}x^{-1}
,\;
bxb^{-1}y^{-1}
>.
$$
Both of these complex
are covered by the Cartesian product of two trees
(Cayley graphs of the free group~$F_2$);
and no finite common cover exists, because
the fundamental
group of such hypothetical covering complex
would embed in both groups $\Gamma_i$
as finite-index subgroups,
but,
in $\Gamma_1$, any finite-index subgroup
contains a finite-index subgroup
which is the
direct product of free groups,
while $\Gamma_2$ has no such
finite-index subgroups [JaW09]
($\Gamma_2$ is not even residually finite
[CaW18],
[BoK21]).
The results of [JaW09] implies also a
minimality of this example in the sense that
\disp{\sl\hfuzz6.7pt
if we restrict ourselves to complexes $K_i$ covered
by products of trees,
then four two-dimensional cells is
the minimum among all non-Leighton pairs.
}%
If we do not restrict ourselves, then
smaller non-Leighton pairs arise.

\proclaim Main theorem \rm(a simplified version).
There exist two finite two-complexes
containing two 2-cells each
that have a common
covering, but have not finite common coverings.

\noindent
(Explicit forms of these two two-2-cell two-complexes
can be found at the very end of this paper.)
Thus, the only question remaining open concerns
complexes with a single 2-cell.
This question seems to be difficult
(although it is closely related to
the well-developed
theory one-relator groups).
The point is that a classification
of one-relator groups up to commensurability is
not an easy task
even for the \emph{Baumslag--Solitar
groups} $\BS(n,m)\:=\pres<c,d|c^{nd}=c^m>$ (though, in this special case,
it was recently obtained [CKZ19]).  Henceforth, $x^{ky}\:=y^{-1}x^ky$,
where $x$ and $y$ are elements of a group and $k\in\Z$.

In conclusion, note that
results on coverings of two-complexes can imply
nontrivial facts about graphs, because one can
``model" 2-cells in graphs by means of additional
vertices and edges, see [BrS21].
Higher dimensional complexes are of little sense here:
if complexes $K_1$ and $K_2$ form a non-Leighton pair,
then their two-skeleta also form such a pair,
as is easy to verify.
A detailed exposition of the general theory of
coverings and CW-complexes
can be found, e.g., in~[FoF89].

\s 1.
Algebraic lemmata

The following fact is well and long known [Mes72],
we give a short proof for the reader's convenience.

\proclaim Commutator lemma.
In the group
$
H=\BS(3,5)=\pres<c,d|c^{3d}=c^5>,
$
the commutator~$h=[c^d,c]$
belongs to any finite-index subgroup.

\Proof
Each finite-index subgroup contains a normal
finite-index subgroup (see, e.g., [KaM82])
Therefore, it suffices to show that $h$ lies in the
kernel of any homomorphism $\phi\:H\to K$ to any finite group $K$.

The elements
$\phi(c^3)$ and $\phi(c^5)$ have the same order
(because they are conjugate); hence, the
order of~$\phi(c)$ is not divisible by three.
Therefore, $\phi(c)\in\gp{\phi(c^3)}$.
Thus,
$\phi(c)^{\phi(d)}\in\gp{\phi(c)}$ and
$h=[c^d,c]$ belongs to the kernel of $\phi$.
This completes the proof.

\proclaim Bottle lemma.
If a group $G$ has a
subgroup
$
\gp{a,b}=\pres<a,b|a^b=a^{-1}>\iso\BS(1,-1),
$
and the element~$b$ lies in all finite-index subgroups of $G$,
then any finite-index subgroup
of $G$
contains a subgroup
isomorphic to the Klein-bottle group
$\BS(1,-1)$.

\Proof
Any finite-index subgroup
contains all elements
conjugate to~$b$,
because the intersection $R$ of
all
finite-index subgroups is normal.
Therefore, $a^2=b^{-1}b^a\in R$ and $\gp{a^2,b}\subseteq R$
It remains to note that
$
a^{2b}
=a^{-2},
$
and the groups~$\gp{a^2}$~and~$\gp b$ are infinite;
hence, the subgroup $\gp{a^2,b}$ is
isomorphic to $\BS(1,-1)$,
because,
\dispno{\sl\hfuzz14pt
in any group, infinite-order elements
$x$ and $y$ such that
$x^y=x^{-1}$ generate a subgroup isomorphic
to the Klein-bottle group.
}(1)%
Indeed, there is obvious epimorphism
$\phi\:\BS(1,-1)=\gp{a,b}\to\gp{x,y}$.
Any element $g\in\BS(1,-1)$
can be written in the form $g=a^kb^l$.
If $g=a^kb^l\in\ker\phi$,
then
$
\ker\phi\ni[b,g]=b^{-1}b^{-l}a^{-k}ba^kb^l=a^{\pm2k}.
$
Therefore, $k=0$ (because $|\gp x|=\infty$).
But then $l=0$ too, because $1=\phi(g)=\phi(b^l)=y^l$,
and~$|\gp y|=\infty$. Thus, $\ker\phi=\1$ and
$\phi$ is an isomorphism.
This completes the proof.

\proclaim No-bottle lemma.
The amalgamated free product
$$
G=
\pres<a,c,d|[a,[c^d,c]]=1,\; c^{3d}=c^5>
=
\pres<a,b|[a,b]=1>\mathop*_{b=[c^d,c]}\pres<c,d|c^{3d}=c^5>
$$
of the free abelian group and the Baumslag--Solitar group
$\BS(3,5)$ contains no subgroups isomorphic to the
Klein-bottle group $K=\BS(1,-1)$.

\Proof
The group $\BS(3,5)$ does not contain
subgroups isomorphic to $K$ [Lev15]
and is torsion-free.
Therefore,
applying once again (1),
we obtain that the
quotient
$
G/\nc{[a,G]}=\gp a_\infty\times\BS(3,5)
$
by the normal closure $\nc{[a,G]}$
of the set~$[a,G]$ of commutators of $a$
and all elements of $G$
has
no nonidentity elements
conjugate their inverse.
Therefore, any element of~$G$
conjugate to its inverse lies in $N=\nc{[a,G]}$.
This subgroup intersects trivially the free factors (and their
conjugates). Thus,
all elements of $N$ have length at least two, and
the following conjugation criterion
(see, e.g., [LS80]) applies:
\disp{\sl
two cyclically reduced words of length $\ge2$
in an
amalgamated free product~$U\mathop*\limits_WV$
are conjugate if and only if one of them can be
obtained from the other by a cyclic permutation
and subsequent conjugation by an element of $W$.
}%
In $U\mathop*\limits_WV$, an equality of reduced words
$u_1v_1\dots=u_1'v_1'\dots$ implies the equalities of the double
cosets $Wu_1W=Wu_1'W$, $Wv_1W=Wv_1'W$, \dots\
Therefore, if a cyclically reduced
word
$
x\in
N\nin\,\pres<a,b|[a,b]=1>\mathop*\limits_{b=[c^d,c]}\pres<c,d|c^{3d}=c^5>
$
is conjugate to its inverse, then, for a letter $x_1$ of $x$, we
obtain the equality $x_1=b^kx_1^{-1}b^l$
(because the map $f_{k,n}\:i\mapsto k-i\pmod n$ from the set
(of subscripts) $\{1,\dots,n\}$ to itself has either a fixed point,
or
an almost fixed point:
$f_{k,n}(i)=i+1\pmod n$ for some $i$;
the latter case would imply that $x_1=b^kx_2^{-1}b^l$
for some adjacent letters $x_1$ and $x_2$ of the reduced
word $x$,
which is impossible).
Substituting $x_1=b^k\^x_1$,
we obtain
$\^x_1^2=b^{l-k}$; thus
\-
either
$\^x_1^2\in\gp{b}$ for some
$\^x_1\in\(\gp{a}_\infty\times\gp{b}_\infty\)\setminus\gp b$,
\smallskip
\-
or
$\^x_1^2\in\gp{[c^d,c]}$ for some
$\^x_1\in\pres<c,d|c^{3d}=c^5>\setminus\gp{[c^d,c]}$.

\enditem
The first is impossible of course. The impossibility of the second case
can
be verified, e.g., as follows.
\-
The quotient
group~$Q=\pres<c,d|c^{3d}=c^5>/\nc{[c^d,c]}$ is torsion-free;
\itemitem{}
indeed, $Q$ is the HNN-extension
$Q=\pres<c,e,d|[e,c]=1,\; e^3=c^5,\; c^d=e>$
of the abelian group
\newline
$A=\pres<c,e|[e,c]=1,\; e^3=c^5>$,
which is torsion-free
(moreover, it is easy to verify that
$A\iso\Z$ and $Q\iso\BS(3,5)$);
\-
therefore,
$\^x_1$ lies in the normal closure $F=\nc{[c^d,c]}$, which is a free
group, because, by the Karrass--Solitar theorem (see, e.g., [LS80]), any
subgroup of an HNN-extension is free if it intersects conjugates of the
base trivially.  It remains to show that $[c^d,c]$ is not a square in $F$
(because in a free group
an inclusion $\alpha^2\in\gp\beta$
implies that $\gp{\alpha,\beta}$ is cyclic
by the Nielsen--Schreier theorem and, hence,
$\alpha\in\gp\beta$ if $\beta$ is not a square).
The commutator~$[c^d,c]$ is not a square in $F$,
because,
assuming the contrary and noting that automorphic images of squares are
squares too, we obtain
$
F=\nc{[c^d,c]}=\nc{\^x^2}\subseteq\gp{\{f^2\;|\;f\in F\}},
$
which cannot hold in a nontrivial free group $F$. This completes the
proof.

\s 2.
Proof of the main theorem

Take the fundamental groups of the torus and the Klein bottle:
$$
G_1=\BS(1,1)=\pres<a,b|[a,b]=1>
\qqbox{and}
G_{-1}=\BS(1,-1)=\pres<a,b|a^b=a^{-1}>
$$
and consider the amalgamated free products
$H_\epsilon=G_\epsilon\*_{b=h}H$
of $G_\epsilon$ and a group $H=\pres<X|R>\supseteq\gp h_\infty$
(henceforth $\epsilon=\pm1$).
Let $K_\epsilon$ be the standard complex of the
(standard)
presentation of $H_\epsilon$:
$$
H_\epsilon=\pres<\{a\}\sqcup X|\{a^{\^h}a^{-\epsilon}\}\sqcup R>,
\qbox{where $\^h$ is a word in the alphabet $X^{\pm1}$
representing the element $h\in H$}.
$$
The Cayley graphs of $G_\epsilon$ are isomorphic surely
(as abstract undirected graphs),
the same is true for the universal coverings of the standard
complexes of presentations of the groups $G_\epsilon$
(these covering complexes are planes partitioned on squares,
Fig. 1).


\def\caption{\vbox{\small
\centerline{%
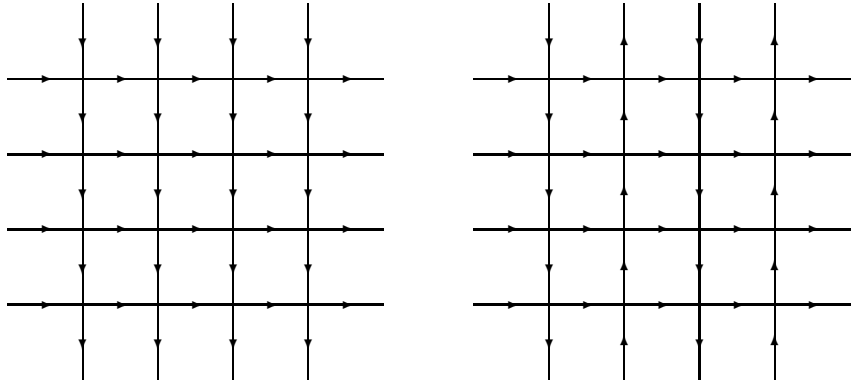
Universal coverings of the standard complexes
of presentations
$G_1$ (left) and $G_{-1}$ (right);
}
\centerline{%
vertical/horizontal edge are labeled by $a$ and $b$,
respectively;
}
\centerline{%
each small square is filled with a 2-cell.
}
}}

\goodbreak
\bigskip
\centerline{\input 1.PIC}
\nobreak%
\centerline{Fig. \lowercase{1}}%
\goodbreak
\bigskip

\noindent
A slightly less trivial observation is that,
for groups $H_\epsilon$,
the universal coverings are isomorphic too:
\dispno{\sl\narrower\narrower
for any infinite-order element $h$ of any group $H$,
the universal coverings of
complexes $K_\epsilon$ are isomorphic.
}(*)%
In what follows, we explain this simple fact in details;
the readers who regard this fact as
obvious, can skip to Observation $(**)$.

It suffices to show that some coverings
$\^K_\epsilon\to K_\epsilon$ have isomorphic $\^K_\epsilon$;
we prefer to take the coverings corresponding to
the normal closure $\nc a$ of $a\in H_\epsilon$.
In explicit form, these complexes $\^K_\epsilon$ are the following ones:
\-
the vertices are elements of $H$;
\-
the edges with labels
from $X$ are drawn as in the Cayley graph of the group~$H$:
an edge with
label $x\in X$
go from each vertex $h'\in H$ to the vertex $h'x\in H$;
\-
in addition, to each vertex $h'\in H$, a directed loop
(edge) $a_{h'}$ labelled by $a$ is attached;
\-
to each cycle whose label is
a relator from $R$, an oriented 2-cell is attached;
\-
to each cycle with label $a^{\^h}a^{-\epsilon}$,
an oriented 2-cell
(a \emph{special} cell) is attached;
thus, going along the boundary of a special cell
in the positive
direction, we meet
two edges labelled by $a$,
namely, $a_{h'}$ and $a_{h'h}^{-\epsilon}$,
where, as usual, $a_{h'h}^{-1}$
means that
the
edge $a_{h'h}$ is traversed against its direction.

\enditem
The isomorphism $\Phi\:{\^K_1}\to\^K_{-1}$
is the following:
\-
the vertices, edges with labels from $X$ and nonspecial 2-cells
(corresponding to relators from $R$) are mapped identically;
\-
to define the mapping $\Phi$ on edges labelled by $a$
and special 2-cells, we choose a set~$T$ of
left-coset representatives
of $\gp h$ in $H$,
and put  $\Phi(a_{th^k})=a_{th^k}^{(-1)^k}$
for all $t\in T$ and $k\in\Z$
(i.e., in each coset, each second
loop labelled by $a$ is inverted);
then the mapping of singular cells are defined naturally:
a cell of $\^K_1$ with
edges~$a_{h'}$~and~$a_{h'h}^{-1}$ on its boundary
is mapped to the cell of $\^K_{-1}$ containing $a_{h'}$ and
$a_{h'h}$ on its boundary.

\medskip

\enditem
The next simple observation is that,

\dispno{\sl
if $h\in H$
belongs to all finite-index subgroups
of $H$, and the complexes $K_\epsilon$ have a finite
common covering,
then the group $H_1$ contains a subgroup isomorphic
to the Klein-bottle group
$\BS(1,-1)$.
}(**)%
Indeed, in $H_{-1}$, the element $b=h$
is contained in all subgroups of finite index
(because the intersection of each such subgroup with $H$
is of finite index in $H$ and, therefore, contains $h$).
By the bottle lemma (applied to $G=H_{-1}$),
we obtain that each finite-index subgroup
contains a subgroup isomorphic to the Klein-bottle group.
It remains to note that,
if a finite complex $\^K$ covers
$K_1$ and $K_{-1}$, then its fundamental group
$\pi_1(\^K)$
embeds into $\pi_1(K_\epsilon)=H_\epsilon$ as a finite-index subgroup.

\medskip

Now, we take a particular group $H$,
namely, let $H$ be the Baumslag--Solitar group:
$
H=\BS(3,5)=\pres<c,d|c^{3d}=c^5>,
$
and let $h\in H$ be the commutator:
$
h=[c^d,c].
$
This element $h$ is contained in any finite-index subgroup
of $H$ by the commutator lemma.
According to $(**)$, this means
that, if complexes~$K_\epsilon$ would have a common finite
covering, then
$H_1=\pres<a,c,d|[a,[c^d,c]]=1,c^{3d}=c^5>$
would contain the Klein-bottle group as a
subgroup,
which contradicts the no-bottle lemma.
Therefore, there are no finite common coverings for complexes~$K_\epsilon$;
while an infinite common covering exists according to $(*)$.
Thus, the following fact
is proven.

\proclaim Main theorem.
The standard complexes of presentations
$
H_\epsilon=\pres<a,c,d|a^{[c^d,c]}=a^\epsilon,\; c^{3d}=c^5>,
$
where $\epsilon=\pm1$,
containing
two 2-cells
\(and one vertex, and three edges\)
have a common covering, but have no finite common coverings.

\References

[AFS91]
J. Abello, M. R. Fellows, J. C. Stillwell,
On the complexity and combinatorics of covering finite complexes,
Australasian Journal of Combinatorics, 4 (1991), 103-112.

[BaK90]
H. Bass, R. Kulkarni,
Uniform tree lattices,
J. Amer. Math. Soc., 3:4 (1990), 843-902.

[BoK21]
I. Bondarenko, B. Kivva,
Automaton groups and complete square complexes,
Groups, Geometry, and Dynamics, 16:1 (2022), 305-332.
\arXiv 1707.00215

[BrS21]
M. Bridson, S. Shepherd,
Leighton's theorem: extensions, limitations, and quasitrees,
Algebraic and Geometric Topology (to appear).
\arXiv 2009.04305.

[CaW18]
P.-E. Caprace, P. Wesolek,
Indicability, residual finiteness,
and simple subquotients of groups acting on trees,
Geometry and Topology, 22:7 (2018), 4163-4204.
\arXiv 1708.04590

[CKZ19]
M. Casals-Ruiz, I. Kazachkov, A. Zakharov,
Commensurability of Baumslag--Solitar groups,
Indiana Univ. Math. J., 70:6 (2021), 2527-2555.
\arXiv:1910.02117

[FoF89]
A. Fomenko, D. Fuchs,
Homotopical topology, 2nd ed.,
Graduate Texts in Mathematics, Vol. 273, Springer,
Cham, 2016.

[JaW09]
D. Janzen, D. T. Wise,
A smallest irreducible lattice in the product of trees,
Algebraic and Geometric Topology, 9:4 (2009), 2191-2201.

[KaM82]
M. I. Kargapolov, Yu. I. Merzljakov,
Fundamentals of the theory of groups,
Graduate Texts in Mathematics, 62, Springer, 1979.

[Lei82]
F. T. Leighton,
Finite common coverings of graphs,
J. Combin. Theory, Series B, 33:3 (1982), 231-238.

[Lev15]
G. Levitt,
Quotients and subgroups of Baumslag--Solitar groups,
J. Group Theory, 18:1 (2015), 1-43.
\arXiv 1308.5122

[LS80]
R. Lyndon, P. Schupp,
Combinatorial group theory,
Springer, 2015.

[Mes72]
S. Meskin,
Nonresidually finite one-relator groups,
Trans. Amer. Math. Soc. 164 (1972), 105-114.

[Neu10]
W. D. Neumann,
On Leighton's graph covering theorem,
Groups, Geometry, and Dynamics, 4:4 (2010), 863-872.
\arXiv 0906.2496

[SGW19]
S. Shepherd, G. Gardam,  D. J. Woodhouse,
Two generalisations of Leighton's Theorem,
arXiv:1908.00830.

[Tuc90]
T. W. Tucker,
Some topological graph theory for topologists:
A sampler of covering space constructions.
In: Latiolais P. (eds)
Topology and Combinatorial Group Theory.
Lecture Notes in Mathematics,
1440 (1990). Springer, Berlin, Heidelberg.

[Wis96]
D. T. Wise,
Non-positively curved squared complexes:
Aperiodic tilings and non-residually finite groups.
PhD Thesis, Princeton University, 1996.

[Wis07]
D. T. Wise,
Complete square complexes,
Commentarii Mathematici Helvetici, 82:4 (2007), 683-724.

[Woo21]
D. Woodhouse,
Revisiting Leighton's theorem with the Haar measure,
Mathematical Proceedings of the Cambridge Philosophical Society,
170:3 (2021), 615-623.
\arXiv 1806.08196

\end